\documentclass[a4paper,10pt]{article}
\usepackage{color,amsmath}
\usepackage{graphics,float}
\usepackage[latin1]{inputenc} 
\usepackage[T1]{fontenc} 
\usepackage{lmodern} 
\usepackage[left=3.2cm,top=3.5cm,right=3.2cm,nohead,foot=1cm]{geometry}
\usepackage[english]{babel} 
\usepackage{tikz}
\usetikzlibrary{arrows,shadows}
\usepackage[affil-it,auth-sc]{authblk}

\usepackage[pdftex]{hyperref} 

\usepackage{amsfonts,latexsym,amssymb} 
\usepackage{mathrsfs,amsmath}
\usepackage[all]{xy}
\graphicspath{{figures/}}

\newcommand{\vecA}{\mathbf{A}}
\newcommand{\vecB}{\mathbf{B}}

\newcommand{\vecF}{\mathbf{F}}

\newcommand{\vecT}{\mathbf{T}}

\newcommand{\vece}{\mathbf{e}}
\newcommand{\vecf}{\mathbf{f}}
\newcommand{\vech}{\mathbf{h}}

\newcommand{\vecg}{\mathbf{g}}

\newcommand{\vecp}{\mathbf{p}}

\newcommand{\vecu}{\mathbf{u}}
\newcommand{\vecv}{\mathbf{v}}

\newcommand{\vecx}{\mathbf{x}}

\newcommand{\vecz}{\mathbf{z}}

\title{Optimal Design for Purcell Three-link Swimmer}
\author{Laetitia Giraldi\footnote{UMA, ENSTA, France laetitia.giraldi@ensta-paristech.fr.edu\vspace{0.1cm}}, Pierre Martinon\footnote{CMAP, Ecole Polytechnique, France}, Marta Zoppello\footnote{University di padova, Italia}}

\begin{document}

\maketitle

\begin{abstract}

In this paper we address the question of the optimal design for the Purcell $3$-link swimmer. 
More precisely we investigate the best link length ratio which maximizes its displacement. 
The dynamics of the swimmer is expressed as an ODE, using the Resistive Force Theory \cite{GrayHancock55}.
Among a set of optimal strategies of deformation (strokes), we provide an  asymptotic estimate of the displacement  for small deformations, 
from which we derive the optimal link ratio. Numerical simulations are in good agreement with this theoretical estimate, 
and also cover larger amplitudes of deformation. Compared with the classical design of the Purcell swimmer,
we observe a gain in displacement of roughly 60\%.

\end{abstract}



\maketitle

\section{Introduction}
The study of self-propulsion at microscopic scale is attracting increasing attention in the recent literature both because of its 
intrinsic biological interest, and for the possible implications on the design of bio-inspired artificial replicas reproducing the 
functionalities of biological systems (see for instance \cite{Lighthill75,BrennenWinet77,PowersLauga09,Bibette}).
At this scale, inertia forces are negligible compared to the viscous ones i.e. low Reynolds number, calling for different swimming 
strategies than at greater scales. Thus, we assume that the surrounding fluid is governed by Stokes equations which 
implies that hydrodynamic forces and torques are linear with respect to the swimmer's velocity. In the case of planar flagellar propulsion
, the Resistive Force Theory (RFT) provides a simple and concise way to compute a local approximation of hydrodynamic forces and Newton 
laws (see \cite{GrayHancock55}). The resulting equations can be written as a system of linear ODEs (see \cite{AlougesDeSimone13,
BermanKenneth13,GMZ13}). In this paper we focus on one of the first example of micro-swimmer model found in literature: 
the ``three-link swimmer'' \cite{Purcell77}. This model is still attracting interest in recent studies, see \cite{BeckerKoehler03,Avron08}.
The structure of the equations of motion leads to establish a connexion between geometrical control theory and micro-swimming (see 
\cite{Montgomery02}). In this paper, we address the optimal design issue, namely finding the optimal length ratio between the three 
links which maximizes displacement of the swimmer. A similar issue has been studied in \cite{TamHosoi07} where a Fourier expansion is 
used to derive an optimal design. Here, techniques from the control theory are used to approximate the leading order term of the 
swimmer's displacement. Maximizing this leading term gives a theoretical value for the optimal link ratio.
As far as we know, this procedure is original in that context, and could be applied to others models such as the three-sphere swimmer 
(see \cite{NajafiGolestanian04}).

The paper is organized as follows.
Section \ref{Sec:Model} recalls the equations of motion for the Purcell swimmer.
Section \ref{sec:Optimalstrategies} presents strokes which maximize the $x$-displacement, based on previous simulations from 
\cite{GMZ13}. Section \ref{sec:optimaldesign} details the expansion of the displacement for such strokes at small amplitude.
By maximizing the leader term of this expansion, we derive an optimal length ratio.
Section \ref{sec:numerical} shows the numerical simulations whose results are consistent with this theoretical ratio, for both small 
and large amplitude of deformation.

\section{Modeling}
\label{Sec:Model}
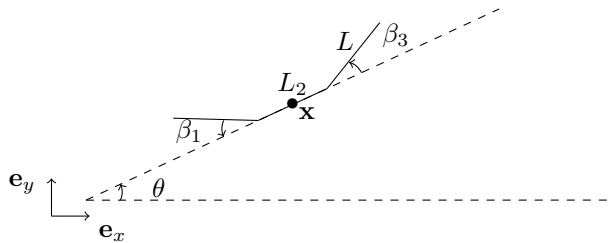
\begin{figure}[h]
\begin{center}
\begin{tikzpicture}[scale=0.5,rotate=25]
\draw[dashed] (-6,0)--(6,0);
\draw (-1,0) -- (1,0);
\draw (-1,0) -- ( -3,1);
\draw (1,0) -- ( 3,1);
\draw [dashed] (-6,0) -- (6.68,-5.91);
\draw [->](-7,0) -- (-7 +0.906,-0.422);
\draw [->] (-7,0) -- (-7 +0.422, 0.906);
\draw (0,0) node {$\bullet$};
\draw [->](2,0) arc (0:50:0.5);
\draw [<-] (-2,0) arc (0:-50:-0.5);
\draw [<-] (-5,0) arc (0:-50:0.5);
\draw (0,0.2) node[below right] {$\vecx$};
\draw (-3,-0.) node[above] {$\beta_1$};
\draw (3,-0.) node[above] {$\beta_3$};
\draw (-4.3,0.1) node[below right] {$\theta$};
\draw (-7 +0.906,-0.422) node[below right] {$\vece_x$};
\draw (-7 +0.422, 1.3)node[below left] {$\vece_y$};
\draw (0,0) node[above] {$L_2$};
\draw (2.2, 0.3)node[above left] {$L$};
\end{tikzpicture}
\caption{Purcell's 3-link swimmer.}
\label{3_links_swimmer}
\end{center}
\end{figure}

\textbf{Purcell's 3-link swimmer.} The $3$-link swimmer is modeled by the position of the center of the second stick $\vecx=(x, y)$, the angle $\theta$ between 
the x-axis and the second stick (the orientation of the swimmer). The shape of the swimmer defined by the two relative angles $\beta_1$ and $\beta_3$ 
(see Fig~\ref{3_links_swimmer}). We also denote by $L$ and $L_2$ the length of the two external arms and central link.

\noindent \textbf{Dynamics via Resistive Force Theory.}
We approximate the non local
hydrodynamic forces exerted by the fluid on the swimmer with local drag forces depending
linearly on the velocity.
We denote by $\vece_i^{\parallel}$ and $\vece_i^{\bot}$ the unit vectors parallel
and perpendicular to the $i$-th link, and we also introduce  $\vecv_i(s)$ the velocity of the point at distance $s$ 
from the extremity of the $i$-th link, that is
\begin{equation*}
\begin{aligned}
&\vecv_1(s)=\dot{\vecx}-\frac{L_2}{2}\dot\theta\vece_2^{\bot}-s(\dot\theta-\dot\beta_1)\vece_1^{\bot},\quad s\in[0,L],\\
&\vecv_2(s)=\dot{\vecx}-(s-\frac{L_2}{2})\dot\theta\vece_2^{\bot},\quad s\in[0,L_2],\\
&\vecv_3(s)=\dot{\vecx}+\frac{L_2}{2}\dot\theta\vece_2^{\bot}+s(\dot\theta-\dot\beta_3)\vece_3^{\bot},\quad s\in[0,L].
\end{aligned}
\end{equation*}

The force $\vecf_{i}$ acting on the $i$-th segment is taken as
\begin{equation}
\vecf_i(s) :=-\xi \left( \vecv_i(s) \cdot \vece_i^\parallel \right) \vece_i^\parallel- \eta \left(\vecv_i(s)
\cdot \vece_i^{\bot}\right) \vece_i^{\bot},
\label{ForceByResistiveTheory}
\end{equation}
where $\xi$ and $\eta$ are respectively the drag coefficients in the directions of
$\vece_i^{\parallel}$ and $\vece_i^{\bot}$. \\
Neglecting inertia forces, Newton laws are written as
\begin{equation}
\label{Newton_laws}
\left\{
\begin{array}{ll}
\vecF = 0 \,,\\
\vece_z \cdot \vecT_{\vecx} = 0\,,
\end{array}
\right.
\end{equation}
where $\vecF$ is the total force exerted on the swimmer by the fluid,
\begin{equation}
\vecF=\sum_{i=1}^N\int_{0}^{L_{i}}\vecf_{i}(s)\,ds\,,
\end{equation}
 and
$\vecT_{\vecx}$ is the corresponding total torque computed with respect to the
central point $\vecx$,
\begin{equation}
\vecT_{\vecx_1}= \sum_{i=1}^N\int_{0}^{L_{i}} \left(\vecx_i(s)-\vecx_1 \right)  \times \vecf_{i}(s) \,ds\,.
\end{equation}
Since the $\vecf_i(s)$ are linear in $\dot{\vecx},\,\dot\theta,\,\dot\beta_1,\,\dot\beta_3$, the system \eqref{Newton_laws} can be rewritten as
\begin{equation}
\vecA(\vecz)\cdot
\left(\begin{matrix}\dot \vecx\\ \dot \theta\end{matrix}\right) -
\vecB(\vecz)\cdot
\left(\begin{matrix}\dot \beta_1\\ \dot \beta_3\end{matrix}\right)=0,
\end{equation}
where $\vecz(t) := (\beta_1,\beta_3,x, y,\theta)(t)^T$.
The matrix $\vecA$ is known as the "Grand Resistance Matrix" and is invertible (see \cite{AlougesDeSimone13}).
Then the dynamics of the swimmer is finally expressed as an ODE system
\begin{equation}
\label{eq:dynamics}
\dot{\vecz}(t)= f(\vecz, \dot\beta_1,\dot\beta_3) = \vecg_1(\vecz(t))\,\dot{\beta}_1(t)+\vecg_2(\vecz(t))\,\dot{\beta}_3(t)\,,
\end{equation}
where $\begin{pmatrix}\vecg_1\left(\vecz\right) &\vecg_2(\vecz)\end{pmatrix}:= \begin{pmatrix}\mathbb{I}_2\\\vecA^{-1} (\vecz)\vecB (\vecz)\end{pmatrix}$ 
with $\mathbb{I}_2$ the $2\times 2$ identity matrix. The literal expression of the $\vecg_i$ is quite complicated (several pages).

\section{Optimal strokes}
\label{sec:Optimalstrategies}
\textbf{Optimal control problem.}
We are interested in finding a periodic sequence of deformations which maximizes the displacement of the swimmer along the x-axis.
More precisely, we optimize both the link length ratio $L_2/L$ and the deformation of the swimmer over time.
Taking the deformation speed $\dot{\beta}_{1|3}$ as control functions, we obtain the optimal control problem
\begin{displaymath}
(OCP) \left\{
\begin{array}{lr}
\max\ x_2(T) \; \text{ s.t.}\\
\dot \vecz(t) = f(\vecz(t),\dot\beta_1,\dot\beta_3) \quad \forall t \in [0,T]\,,\\
\dot \beta_{1|3}\in \mathbf{U} = [-b,b] \quad \forall t \in [0,T]\,,\\
\beta_{1|3}(t) \in [-a,a] \quad \forall t \in [0,T]\,,\\
x_2(0) = y_2(0) = \theta_2(0) = 0 , y_2(T) = \theta_2(T) = 0\,,\\
\beta_{1|3}(0) = \beta_{1|3}(T),\\
2L + L_2 = c.\\
\end{array}
\right.
\end{displaymath}

We set the constraints $a$ and $b$ over the amplitude and deformation speed, as well as the total length $c$ of the swimmer.
The final time $T$ is fixed, and the constraint $\beta_{1|3}(0) = \beta_{1|3}(T)$ ensures that the swimmer is in the same configuration 
at the initial and final time. Note that this condition can be satisfied by either a single stroke or a sequence of strokes. 
From \cite{GMZ13}, numerically solving $(OCP)$ typically gives a periodic sequence of identical strokes.
Their phase portrait is octagonal, as illustrated on Fig.\ref{Fig:strokes}, and we will detail how this shape is consistent with optimal control theory.\\

\textbf{Pontryagin's Maximum Principle (PMP).}
We recall here the PMP as it gives some insight on the shape of optimal strokes.
This theorem in optimal control introduced by Pontryagin et al. in \cite{Pontryagin} gives necessary conditions for 
local optimality. Interested readers can find more information on the PMP in \cite{Agrachev08,Trelat05}. 
The PMP is characterized by an Hamiltonian function $H$ that formally depends on the state variables 
$\vecz$, the control functions $\dot\beta_{1|3}$, and so-called \emph{costate} variables noted $\vecp$.
While originally inspired by the Hamiltonian in mechanics, in the context of optimal control $H$ does not actually correspond to the 
energy of the system. The costate variables play the part of the generalized velocities in Lagrangian mechanics, and they can be
interpreted as Lagrange multipliers (in the sense of constrained optimization) related to the dynamics of the system.
Let the Hamiltonian be
\begin{equation}
\begin{aligned}
\label{eq:Hamiltonian}
H(\vecz,\vecp,\dot\beta_1,\dot\beta_3)= \left\langle\vecp,\vecg_1(\vecz)\right\rangle \dot\beta_1+\left\langle\vecp,\vecg_2(\vecz)\right\rangle \dot\beta_3.
\end{aligned}
\end{equation}
Under the assumption that $\vecg_{1|2}$ are continuous and $C^1$ with respect to $\vecz$, the PMP states that:\\
if $(\vecz^*,\dot\beta_1^*,\dot\beta_3^*)$ is a solution of $(OCP)$ then there exists $\vecp^* \neq 0$ absolutely continuous such that
$\dot{\vecz}^* =   H_p(\vecz^*,\vecp^*,\dot\beta_1^*,\dot\beta_3^*)$, 
$\dot{\vecp}^* = - H_z(\vecz^*,\vecp^*,\dot\beta_1^*,\dot\beta_3^*)$, 
$\vecp^*(T)$ is orthogonal to the cotangent cone of the final conditions at $\vecz^*(T)$ and
$(\dot\beta_1^*,\dot\beta_3^*)$ maximizes the Hamiltonian for almost every time $t \in [0,T]$.\\

\textbf{Bang arcs.} The Hamiltonian in \eqref{eq:Hamiltonian} is linear in the controls $\dot\beta_{1|3}$.  
If we assume $\langle \vecp,\vecg_i(\vecz)\rangle \neq 0$ for $i=1,2$ over a time interval, then the optimal control $\dot\beta_{1|3^*}$
that maximizes $H$ must be on the boundary of $U =\{(-b,-b),(-b,b),(b,-b),(b,b)\}$. In terms of phase portrait, this corresponds to 
diagonal lines.\\

\textbf{Constrained arcs.} Moreover, we have the constraints on the joint angles $\beta_{1|3}(t) \in [-a,a]$. When one of them is active
and $|\beta_i|= a$, the corresponding control $\dot\beta_i = 0$. In terms of phase portrait, this gives horizontal or vertical lines.\\

{\textbf{Symmetries.} As stated in \cite{TamHosoi07}, we expect optimal strokes to be symmetric with respect to the 
diagonal axes $\beta_1 = \beta_3$ and $\beta_1 = -\beta_3$. This comes from the equations of motion being linear and time independent.
From the linearity, optimal strokes should be invariant by reflection with respect to the axis of the swimmer's body.
From time independence, the stroke should be invariant when inverting the arms movement and going backwards in time.
}

\section{Optimal swimmer design}
\label{sec:optimaldesign}
\begin{figure}[H]
\begin{center}
\includegraphics[width=6cm]{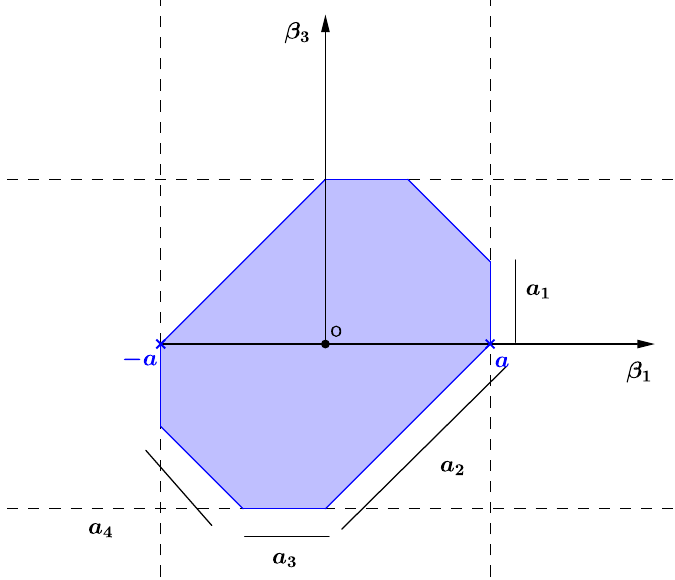}
\caption{\label{Fig:strokes} Phase portrait $(\beta_1,\beta_3)$ of the octagonal stroke considered for the expansion of the displacement.}
\end{center}
\end{figure}

In this section, we express the leader term of the swimmer's displacement for a stroke of small perimeter which satisfies all properties stated in the previous section. We represent the stroke by a closed octagonal curve $\gamma$ in the phase portrait $(\beta_1,\beta_3)$, see Fig.~\ref{Fig:strokes}. 

As a consequence of neglecting inertia forces, velocities appear linearly in the dynamic, and time can be rescaled without changing the dynamics.
Thus the displacement of the swimmer after one stroke does not depend on the speed along the curve $\gamma$, but only on the shape of the stroke.
From now on, we parametrize $\gamma$ by the arc-length $s$.
Using a similar approach to \cite{Coron56}, we express the  swimmer's displacement along the $x$-axis (i.e., $x(T)-x(0)$) as an asymptotic 
expansion for small length $a_i$, $i=1,\cdots,4$.

\textbf{Displacement over the arc $s\in[0,a_1]$.} On this part, according to Fig. \ref{Fig:strokes}, 
we set $\vecu = (\dot \beta_1,\dot \beta_3) = (0,-1)$. The dynamics of the swimmer is
therefore given by $\dot \vecz = - \vecg_2$, and the time expansion
at order two is given by

\begin{eqnarray}
\label{eq:z_at_t1}
\vecz(a_1)=&&\,\vecz(0)-a_1\vecg_2(\vecz(0))\nonumber\\
&&+\frac{a_1^2}{2}
\frac{\partial\vecg_2}{\partial \vecz}{\vecz(0)}\left(\vecg_2(\vecz(0))\right)+ o(a_1^3)\,.
\end{eqnarray}

\textbf{Displacement over the arc $s\in[a_1,a_1+a_2]$.} Similarly, the position of the swimmer at $s= a_1+a_2$ can be expressed as
\begin{eqnarray}
\label{eq:z_at_t1t2}
\vecz(a_1+a_2) = && \vecz(a_1)-\frac{a_2\sqrt{2}}{2}\vech(\vecz(a_1))\nonumber\\
&&+\frac{a_2^2}{4}\frac{\partial\vech}{\partial\vecz}{\vecz(a_1)}\left(\vech(\vecz(a_1))\right)+ o(a_2^3)\,,
\end{eqnarray}
where  $\vech := \vecg_1+\vecg_2$. Plugging the value of $\vecz(a_1)$ from \eqref{eq:z_at_t1} into \eqref{eq:z_at_t1t2} and neglecting the terms of order greater
than two, we get
\begin{eqnarray}
\vecz(a_1+a_2)=&&\vecz(0)+c_1(\vecg_1,\vecg_2,\vecz(0),a_1,a_2)\nonumber\\
&&+c_2(\vecg_1,\vecg_2,\vecz(0),a_1,a_2)\nonumber\\
&&+ o(a_1^3)+o(a_2^3)\,
\end{eqnarray}
with
\begin{eqnarray*}
c_1(\vecf,\vecg,\vecz, a_1,a_2) =&& -\frac{\sqrt{2}a_2}{2} \vecf(\vecz) \\
&&+(-a_1-\frac{\sqrt{2}a_2}{2}) \vecg(\vecz),\\
c_2(\vecf,\vecg,\vecz, a_1,a_2) =&&\frac{a_2^2}{4}\frac{\partial\vecf}{\partial\vecz}{\vecz}\left(\vecf(\vecz)\right)
+\frac{a_2^2}{4}\frac{\partial\vecg}{\partial\vecz}{\vecz}\left(\vecf(\vecz)\right)\\
&&+\left(\frac{a_1a_2\sqrt{2}}{2}+\frac{a_2^2}{4}\right)\frac{\partial\vecf}{\partial\vecz}{\vecz} \left(\vecg(\vecz)\right)  \\
&&+\left(\frac{a_1a_2\sqrt{2}}{2}+\frac{a_2^2}{4}+\frac{a_1^2}{2}\right)\frac{\partial\vecg}{\partial\vecz}{\vecz}\left(\vecg(\vecz)\right).\\
\end{eqnarray*}

\textbf{Displacement over the complete stroke.} Iterating the computations along each arc and noting by $P=2(a_1+a_2+a_3+a_4)$ the stroke perimeter, 
the expansion of the total displacement for the octagonal stroke is finally obtained as
\begin{equation}
\vecz(T) - \vecz(0) = C\ [\vecg_1,\vecg_2](\vecz(0)) + o(a_i^3)_{i=1-4}\,,
\end{equation}
where
$$
C=\frac{a_1a_2\sqrt{2}}{2}+a_1a_3+\frac{a_2a_3\sqrt{2}}{2}+\frac{a_1a_4\sqrt{2}}{2}+a_2a_4+\frac{a_3a_4\sqrt{2}}{2}
$$
and 
$$
[\vecg_1,\vecg_2](\vecz(0)) = \nabla\vecg_2(\vecz(0))\cdot\vecg_1(\vecz(0)) - \nabla\vecg_1(\vecz(0))\cdot\vecg_2(\vecz(0))
$$
is the Lie brackets of $\vecg_1$ and $\vecg_2$ at point $\vecz(0)$.
Choosing the starting point $\vecz(0)$ such that $\theta(0) = \beta_1(0) = \beta_3(0) = 0$, we compute the Lie bracket with a formal calculus tool
\begin{equation}
[\vecg_1,\vecg_2](0,0,x,y,0)=\begin{pmatrix}0\\0\\ \frac{\eta - \xi}
{\xi}\frac{L^3L_2(3L+2L_2)}{(2L+L_2)^4}\\ 0\\ 0\end{pmatrix}\,.
\end{equation}
Consequently, the $x$-displacement after one stroke is approximated by
\begin{equation}
\label{eq:deltax}
x(T) - x(0) = C \left(\frac{\eta - \xi}{\xi}\right) \left(\frac{L^3L_2(3L+2L_2)}{(2L+L_2)^4} \right) + o(a_i^3)_{i=1-4}
\end{equation}
Setting the total length of the swimmer by a constant equal to $c$, i.e., $2L+L_2=c$, we find that
\eqref{eq:deltax} has a unique maximum at
\begin{equation}
L^*=c\Bigl(1-\sqrt{\frac{2}{5}}\Bigr)\,, \quad L_2^*=c\Bigl(2\sqrt{\frac{2}{5}}-1\Bigr)\, ,
\end{equation}
which gives an optimal ratio of
\begin{equation}
\label{eq:formula_ratio}
\left(\frac{L_2}{L}\right)^* = \frac{\sqrt{10}-1}{3}
 \sim 0.721\,.
\end{equation}

\textbf{Remark}: in \cite{TamHosoi07} an optimal ratio of $0.747$ is given for an efficiency-type
criterion. The small gap may be due to the difference in models, or the change of the objective function.

\section{Numerical simulations}
\label{sec:numerical}

We solve now the optimal control problem $(OCP)$ numerically, in order to determine the optimal swimming strategy and link ratio.
Simulations are performed with the toolbox \textsc{Bocop} (\cite{Bocop14}) that implements a direct transcription method.
This approach uses a time discretization to transform the continuous $(OCP)$ into a finite-dimensional optimization problem
(nonlinear programming).
We refer interested readers to \cite{Betts01} for more details on these methods. 
We use here an implicit midpoint discretization with $100$ to $2500$ time steps.
Note that this method does not use the PMP.\\

As stated in $(OCP)$, the criterion is to maximize the total displacement along the x-axis over a fixed time $T$.
The initial state of the swimmer is set as $x(0)=y(0)=\theta_2(0)=0$, with the final conditions $y(T)=\theta_2(T)=0$.
The initial shape angles are left free, with the periodicity conditions $\beta_i(0)=\beta_i(T),i=1,3$.
We set the total length $c=4$ for an easier comparison with the classical Purcell swimmer ($L=1,L_2=2)$.\\

We explore different values for the bounds $a,b$ on the shape angles and deformation speed and see their influence on the optimal stroke
and link ratio.
For practical applications, the values for $a$ and $b$ should reflect the physical characteristics of the studied swimmer.
It should be pointed out that the period of the optimal stroke is not known a priori.
We arbitrarily set $T=1$ in the first set of simulations, and $T=25$ when studying the larger amplitudes.
In the latter case we find that the swimming strategy consists in a periodic sequence of identical strokes, as previously observed in 
\cite{GMZ13}.

\subsection{Small amplitudes, influence of speed limits}

\begin{figure}[H]
\begin{center}
\includegraphics[width=0.7\textwidth]{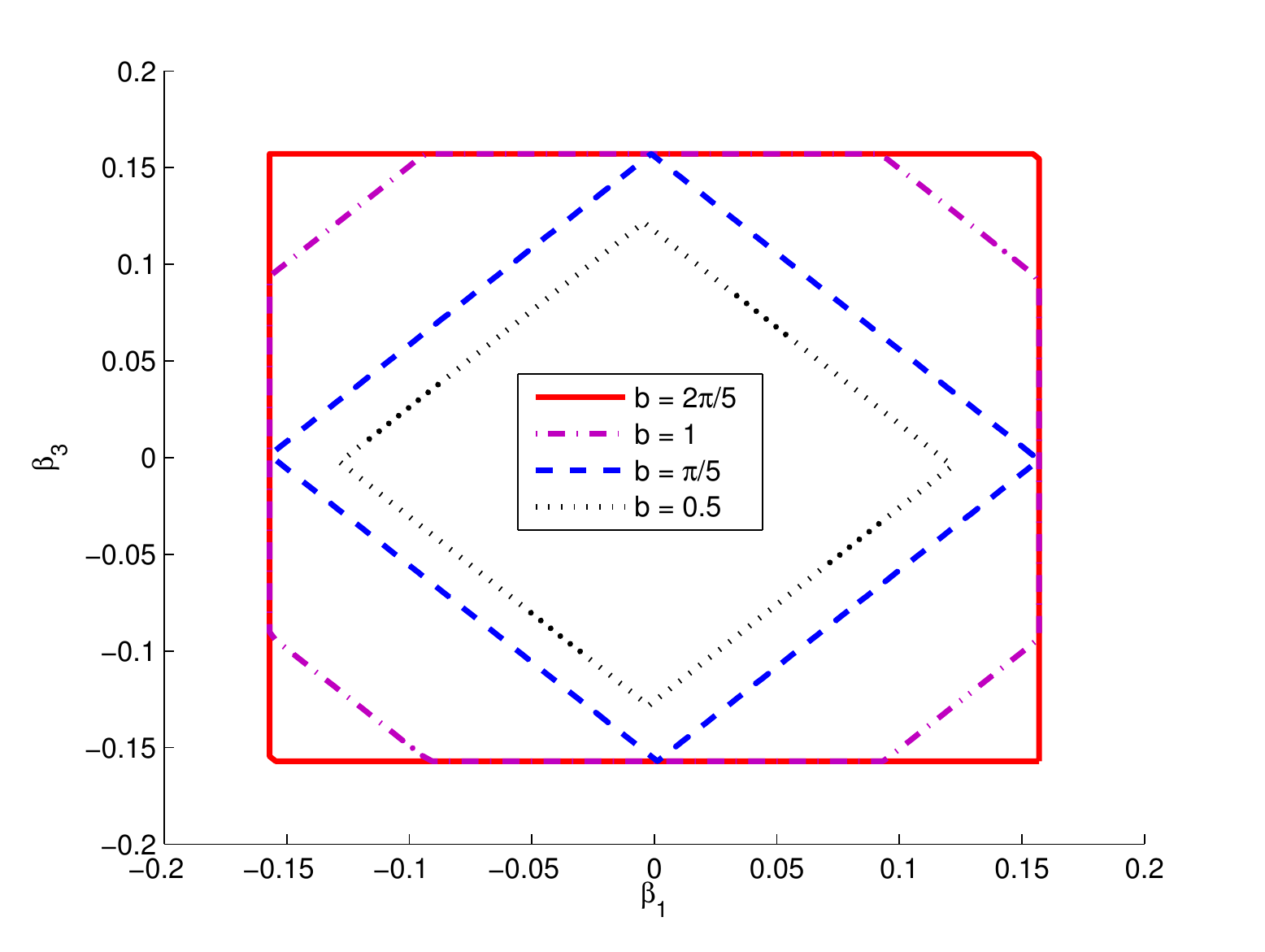}
\caption{\label{fig:portraits-api20} Phase portraits of the strokes for small amplitudes, $a=\pi/20$. 
The shapes observed are consistent with the discussion in section \ref{sec:Optimalstrategies}.}
\end{center}
\end{figure}

We start with small amplitudes by setting $a=\pi/20$ and solve $(OCP)$ for different values of the speed limit $b$.
Here we set $T=1$ and use $250$ time steps for the discretization. Optimizations take about one minute on a standard laptop.
Results are given in Table.\ref{tab:api20}, with the phase portraits for the shape angles $\beta_1,\beta_3$ on 
Fig.\ref{fig:portraits-api20}.\\

First, we observe that the optimal ratio $L_2/L$ is very close to its theoretical value of $0.721$ from 
\eqref{eq:formula_ratio}, regardless of $b$.
The speed bound does however have an influence on the shape of the optimal stroke, and its displacement.
Displacement increases with higher speeds, and we find the following empirical relation between $b$ and the stroke shape, 
confirmed by simulations with other values of $a$:\\
- for $b < 4a/T$: diamond stroke, which touches the bound $a$ for the limit case $b = 4a/T$.\\
- for $4a/T < b < 8a/T$: octagonal stroke.\\
- for $b=8a/T$: classical Purcell stroke (square).\\
- for $b>8a/T$: sequence of several strokes.\\
The three strokes observed (diamond, octagon, square) match the discussion from Section \ref{sec:Optimalstrategies}. They include only diagonal lines 
(bang arcs saturating the speed limit $b$) and horizontal/vertical lines (constrained arcs for the amplitude limit $a$).
Note also that the square and diamond strokes are particular cases of the octagonal one, by setting the appropriate arc lengths to 0.\\

\emph{Remark: this empirical relation can also be interpreted in terms of the period $T$, with the two limit values $T=8a/b$ for the 
Purcell stroke and $T=4a/b$ for the diamond touching $a$.}

\begin{table}[H]
\begin{center}
\caption{\label{tab:api20} Small amplitude ($a=\pi/20$).}
\begin{tabular}{|l|rrr|}
\hline
$b$	&$x(T)$		&$L_2/L$	&stroke\\
\hline
0.5	&2.68E-3	&0.719	&diamond\\
$\pi$/5	&4.23E-3	&0.719	&diamond\\
0.75	&5.70E-3	&0.719	&octagon\\
1	&7.73E-3	&0.719	&octagon\\ 
2$\pi$/5&8.42E-3	&0.717	&square\\
1.5	&1.14E-2	&0.719	&octagon (x2)\\
2	&1.55E-2	&0.719	&octagon (x2)\\
\hline
\end{tabular}
\end{center}
\end{table}

\subsection{Comparison with the classical Purcell swimmer}

\begin{figure}[H]
\begin{center}
\includegraphics[width=0.7\textwidth]{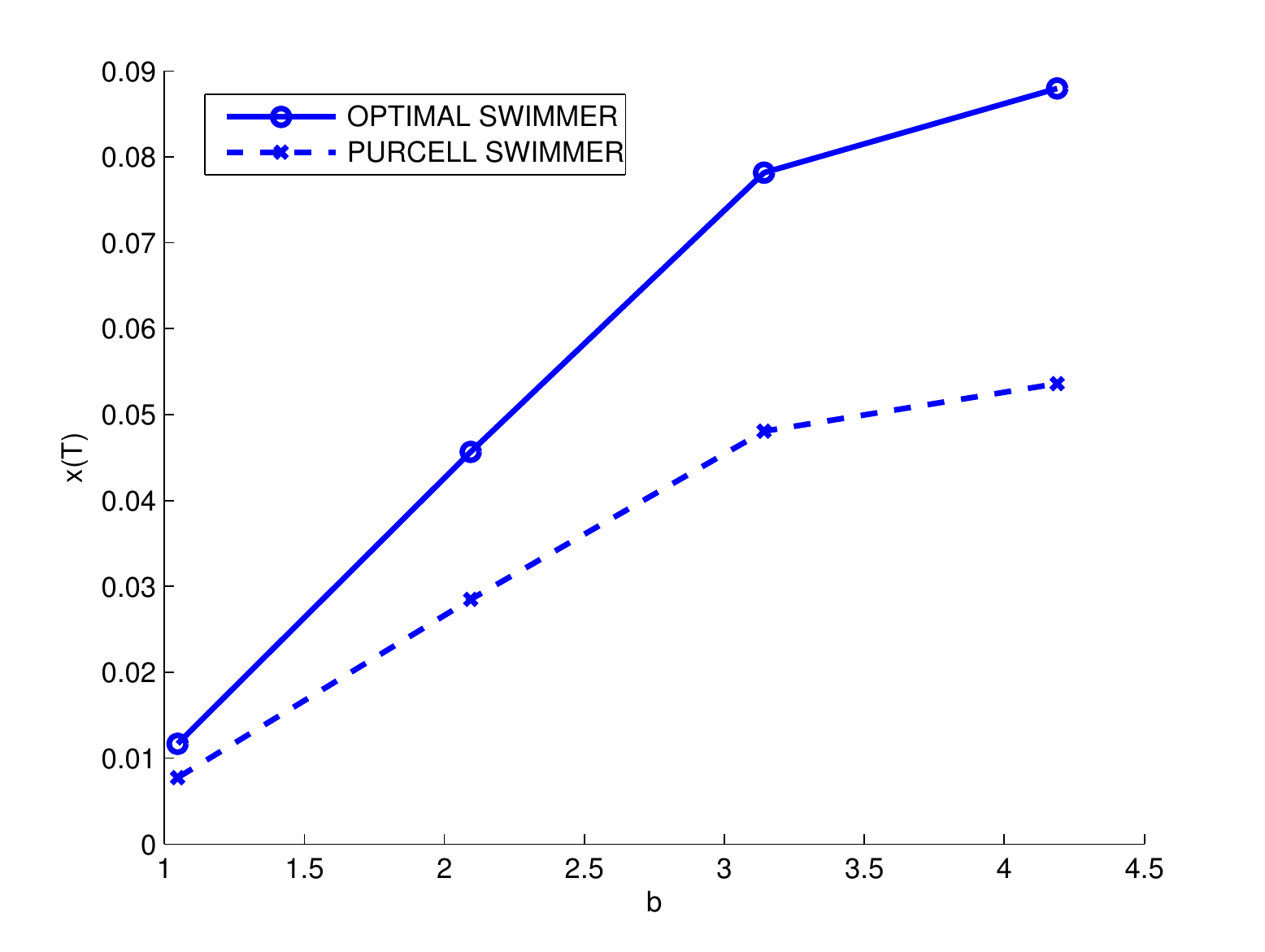}
\caption{\label{fig:OptimalvsPurcell} Displacement for the optimal/Purcell swimmer.}
\end{center}
\end{figure}

Now we compare the performance of the optimal swimmer with respect to the classical Purcell swimmer defined by $L=1,L_2=2$, meaning a ratio of $2$.
For this comparison we set $a=\pi/6$ (thus a stroke amplitude of $\pi/3$) and $b={\pi/3,2\pi/3,\pi,4\pi/3}$ and $T=1$.
The optimization for the Purcell swimmer is done by setting $L=1$ instead of letting it free.
The results are summed up in Table.\ref{tab:api6} and Fig.\ref{fig:OptimalvsPurcell}.
We see that the shape of the stroke matches the empirical law, and that the optimal link ratio stays close to its theoretical value.
We also observe a consistent gain in displacement that seems to increase with the speed limit, up to $64\%$ for the classical Purcell 
stroke (square).

\begin{table}[H]
\begin{center}
\caption{\label{tab:api6} Optimal swimmer vs Purcell swimmer.}
\begin{tabular}{|l|rrr|ll|}
\hline
b	&x(T)	&$L_2/L$	&stroke	&$x_{Purcell}(T)$&gain\\
\hline
$pi/3$	 &1.17E-2 &0.717 &diamond & 7.373E-3 &51\%\\
$2\pi/3$ &4.57E-2 &0.708 &diamond & 2.848E-2 &60\%\\
$\pi$	 &7.82E-2 &0.699 &octagon & 4.806E-2 &63\%\\
$4\pi/3$ &8.80E-2 &0.695 &square  & 5.359E-2 &64\%\\
\hline
\end{tabular}
\end{center}
\end{table}

\subsection{Large amplitudes, influence of angle limits}
Now we study the influence of the maximal amplitude of the stroke, set by the bound $a$.
In this last part we set the deformation speed limit $b=1$ to focus on the amplitude.
Since we would like to study only the true optimal strokes, whose period is not known, we also take a longer final time $T=25$.
We expect to obtain trajectories that exhibit a sequence of several identical strokes with a period $T^* < T$.
The number of time steps is raised accordingly to 2500, which increases the computational time up to half an hour.
Another way of finding the optimal stroke directly could be to leave the final time $T$ free in the optimization, while 
maximizing the average speed of the stroke $x(T)/T$ instead of the displacement $x(T)$.\\

\begin{figure}[H]
\begin{center}
\includegraphics[width=0.7\textwidth]{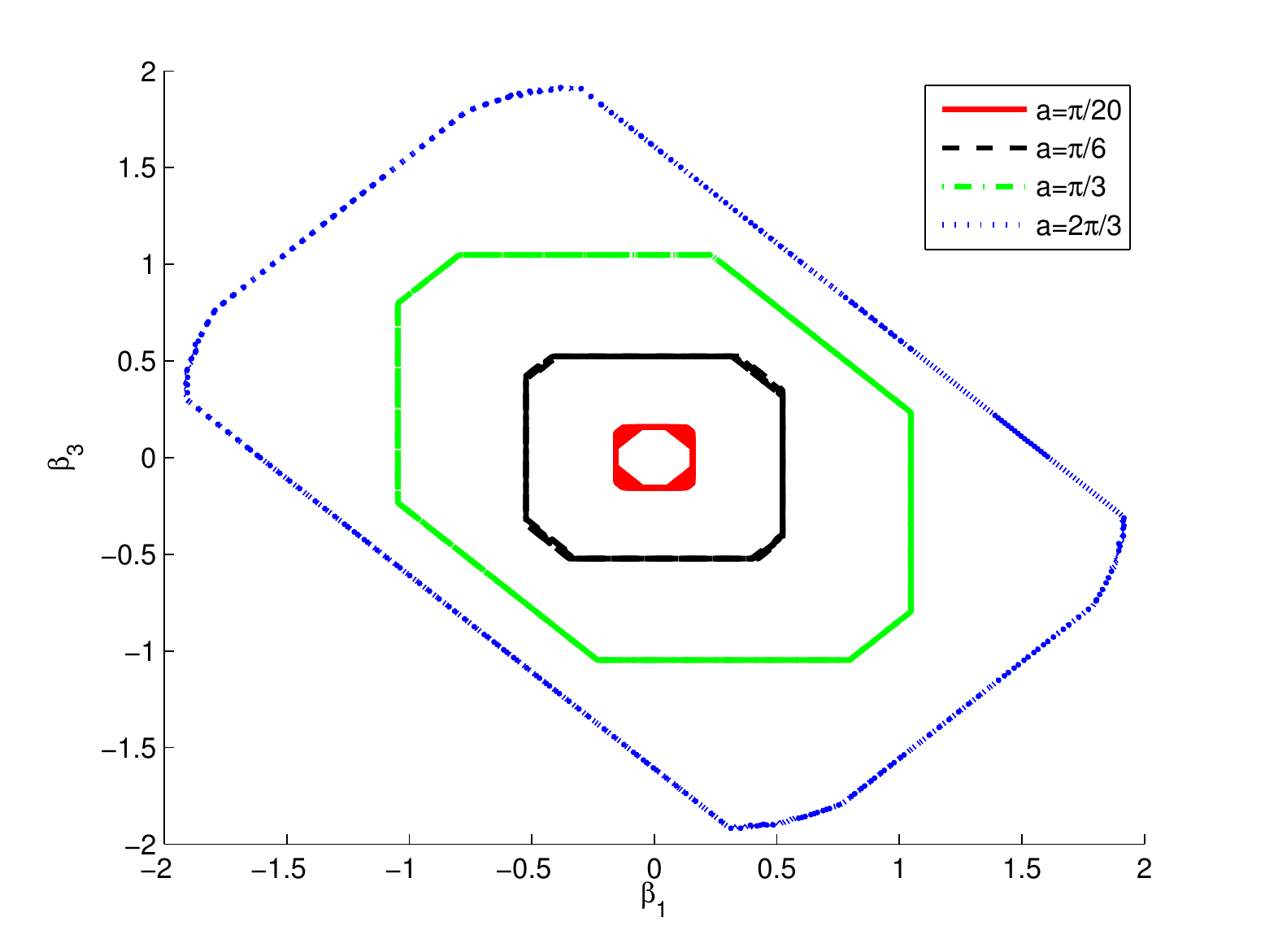}
\caption{\label{fig:portrait-T25}Larger amplitudes - Phase portrait (with several superposed strokes for each trajectory).}
\end{center}
\end{figure}

\begin{figure}[H]
\begin{center}
\includegraphics[width=0.7\textwidth]{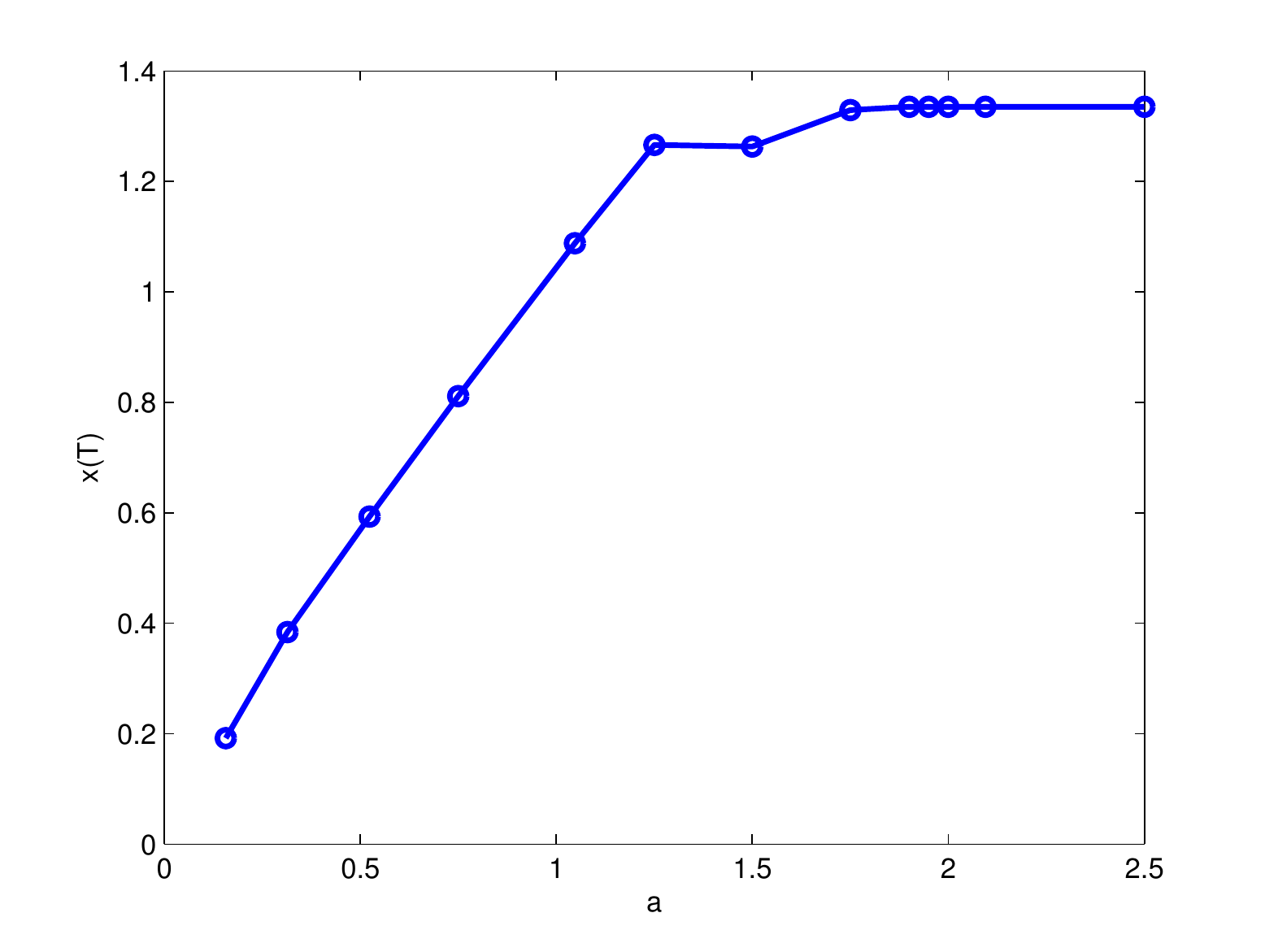}
\caption{\label{fig:displacement-T25} Larger amplitudes - Overall displacement. Note that since the displacement is expected to be a strictly increasing function of $a$, we see that for $a=1.5$, the optimization converged to a local solution.}
\end{center}
\end{figure}

The results are illustrated in Table.\ref{tab:T25} and Figs.\ref{fig:portrait-T25}-\ref{fig:displacement-T25}.
First, the simulations confirm that the optimal strategy is a periodic sequence of identical strokes.
The shape of the optimal stroke is always octagonal until it becomes unconstrained for very large values of $a$.   
We observe that the central symmetry observed for small amplitudes is lost for larger $a$, however symmetry w.r.t both diagonal axes 
still holds as expected.\\

In the unconstrained case, we see arcs that are neither bang arcs (diagonal) or constrained arcs (horizontal/vertical), but rather
appear as smooth curves (see Fig.\ref{fig:portrait-T25}) . These are characteristic of so-called \emph{singular arcs}, namely the case where $\langle p,g_i(z) \rangle=0$
in the PMP. More details on the analysis of singular arcs can be found in \cite{Trelat05}, unfortunately here the complexity of the 
$g_i$ makes further study quite difficult.\\

The total displacement $x(T)$ increases with $a$, first almost linearly when $a<\pi/3$ (see Fig.\ref{fig:displacement-T25}). 
From $a\approx 1.95$ and above, we obtain the same, unconstrained solution.
The improvement in displacement appears to be marginal between $a=\pi/3$ and the unconstrained case.
Note that since the displacement is expected to be a monotone increasing function of $a$, we see that for $a=1.5$, the optimization converged to a local solution.\\

The optimal ratio $L_2/L$ shows a steady decrease with $a$, starting quite close to the value $0.721$ computed for 
small amplitudes, the seemingly reaching a limit value of $2/3$ in the unconstrained case (i.e. $L=1.5,L_2=1$).
We recall that the classical Purcell swimmer has a link ratio of $2$ ($L=1,L_2=2$).

\begin{table}[H]
\begin{center}
\caption{\label{tab:T25} Larger amplitudes: optimal link ratio and stroke. Solutions become unconstrained about $a=1.95$.}
\begin{tabular}{|l|rrr|}
\hline
a		&x(T)	&$L_2/L$	&stroke\\
\hline
$\pi/20$	&0.192	&0.719	&octagon x26\\
$\pi/10$	&0.384	&0.712	&octagon x13\\
$\pi/6$ 	&0.593	&0.697	&octagon x7\\
0.75		&0.811	&0.676	&octagon x5\\
$\pi/3$ 	&1.088	&0.660	&octagon x4\\
1.25		&1.266	&0.660	&octagon x4\\
1.5		&1.263	&0.660	&octagon x3\\ 
1.75		&1.329	&0.667	&octagon x3\\
$2\pi/3$	&1.335	&0.667	&unconstrained x3\\
2.5		&1.335	&0.667	&unconstrained x3\\
\hline
\end{tabular}
\end{center}
\end{table}

\section{Conclusion} 
This study is devoted to the optimization of the link ratio of the three-link swimmer for maximal displacement.
We provide an estimate of the displacement based on an expansion at small deformations, which gives a theoretical optimal link ratio.
Numerical simulations when solving the optimal control numerically are consistent with this theoretical ratio for small amplitudes.
We also observe that the optimal ratio changes for large amplitudes, with a limit value of $0.667$ in the unconstrained case versus
a theoretical ration of $0.721$ at small amplitudes. For an amplitude of $\pi/3$, the displacement gain is about 60\% compared with
the classical Purcell swimmer design.
A possible continuation of this work is the comparison of different objective functions, such as average speed or efficiency.

\bibliographystyle{plain}
\bibliography{Bib_Optimal_Swimmer}

\end{document}